\documentclass[twoside, a4paper, 11pt, article]{memoir}

\usepackage[english]{babel}

\usepackage[latin1]{inputenc}
\usepackage[T1]{fontenc}

\usepackage{microtype}

\usepackage{amsmath}
\usepackage{amsfonts}
\usepackage{amssymb}

\usepackage[amsmath,thmmarks,hyperref]{ntheorem}

\usepackage{graphicx}

\usepackage{mathtools}

\usepackage[ps,all]{xy}
\SelectTips{cm}{}
\CompileMatrices

\def\cxymatrix#1{\xy*[c]\xybox{\xymatrix#1}\endxy}

\usepackage{url}

\usepackage{caption}

\usepackage[justification=RaggedRight]{subfig}

\usepackage[sc]{mathpazo}

\usepackage{varioref}

\ifpdf
\usepackage[pdftitle={The Mapping Class Group Orbit of a Multicurve},
pdfauthor={Rasmus Villemoes},
pdfkeywords={mapping class group, multicurve, group cohomology},
colorlinks=true, 
pdftex, backref=none]{hyperref}
\else
\usepackage[
pdftitle={The Mapping Class Group Orbit of a Multicurve}, 
pdfauthor={Rasmus Villemoes},
pdfkeywords={mapping class group, multicurve, group cohomology},
colorlinks=true, 
ps2pdf,
bookmarksnumbered=true,
backref=none]{hyperref}
\fi
\usepackage{memhfixc}

\newcommand{\SL}{\mathrm{SL}}
\newcommand{\PSL}{\mathrm{PSL}}

\newcommand{\setC}{\mathbb{C}}
\newcommand{\setZ}{\mathbb{Z}}

\DeclareMathOperator{\fut}{fut}
\DeclareMathOperator{\pas}{pas}

\newcommand{\inv}{^{-1}}

\newcommand{\eqvdots}{\mathmakebox[\widthof{${}={}$}][c]{\vdots}}

\theoremstyle{plain}
\theoremsymbol{}
\theorembodyfont{\itshape}
\theoremheaderfont{\normalfont\bfseries}
\theoremseparator{.}
\newtheorem{theorem}{Theorem}[chapter]
\newtheorem{lemma}[theorem]{Lemma}
\newtheorem{proposition}[theorem]{Proposition}

\theoremstyle{nonumberplain}
\theorembodyfont{\normalfont}
\theoremheaderfont{\itshape}
\theoremsymbol{\ensuremath{\square}}
\newtheorem{proof}{Proof}

\qedsymbol{\ensuremath{\square}}

\makepagestyle{orbit}
\makeoddhead{orbit}{}{\textsc{The mapping class group orbit of a multicurve}}{\thepage}
\makeevenhead{orbit}{\thepage}{\textsc{Rasmus Villemoes}}{}
\makeatletter
\makepsmarks{orbit}{%
      \let\@mkboth\markboth
      \def\chaptermark##1{%
        \markboth{{%
          \ifnum \c@secnumdepth >\m@ne
            \if@mainmatter
            \thechapter\ $\cdot$\ %
            \fi
          \fi
          ##1}}{}}%
      \def\sectionmark##1{%
        \markright{{%
          \ifnum \c@secnumdepth > \z@
            \thesection\ $\cdot$\ %
          \fi
          ##1}}}%
      \def\subsectionmark##1{%
        \markright{{%
          \ifnum \c@secnumdepth > \@ne
            \thesubsection\ $\cdot$\ %
          \fi
          ##1}}}%
    }
\makeatother
\title{The mapping class group orbit of a multicurve}
\author{Rasmus Villemoes}
\date{\today}

\begin{document}

\maketitle

\begin{abstract}
  \noindent Given a set equipped with a transitive action of a group,
  we define the notion of an almost invariant coloring of the set. We
  consider the mapping class group orbit of a multicurve on a compact
  surface, and prove that in the case of genus at least two, no such
  almost invariant coloring exists. Conversely, in the case of a
  closed torus, one may find almost invariant colorings using 
  arbitrarily many colors.
\end{abstract}

\chapter{Introduction}
\label{cha:introduction}

Let $G$ be a group and $X$ an infinite set on which $G$ acts. We
define a \emph{coloring} (or $C$-coloring) of $X$ to be any map
$c\colon X\to C$ into some set $C$ of ``colors''. We will use the
following terminology:
\begin{itemize}\tightlist
\item A coloring $c$ is \emph{invariant} if $c(gx) = c(x)$ for
  each $g\in G$ and $x\in X$.
\item A coloring is \emph{almost invariant} if, for each $g\in G$, the
  identity $c(x) = c(gx)$ fails for only finitely many $x\in X$.
\item Two colorings are \emph{equivalent} if they assign different colors to
  only finitely many elements of $X$; this is clearly an equivalence
  relation on the set of $C$-colorings.
\item A coloring is \emph{trivial} if it is equivalent to a
  monochromatic (constant) coloring.
\end{itemize}
We will only deal with the case where the action of $G$ is
transitive. Then clearly the only invariant colorings are the constant
ones, and hence we are only interested in studying the question of
existence of almost invariant colorings. If two colorings are
equivalent and one is almost invariant, so is the other, which
explains the above definition of a trivial coloring. If one wants to
classify all almost invariant colorings, this can clearly not be done
better than up to the equivalence defined above.

A \emph{simplification} of $c$ is a coloring obtained by
post-composing $c$ with some map $i\colon C\to C'$ (one ``identifies''
some of the colors). Clearly a simplification of an almost invariant
coloring is almost invariant. Now, if there exists an almost invariant,
non-trivial $C$-coloring $c$, there also exists an almost invariant
coloring where exactly two colors are used.  To see this, partition
$C$ into $C_0\sqcup C_1$ such that $c\inv(C_k)$, $k=0,1$, are both
infinite, and define a $\{0,1\}$-coloring by composing $c$ with the
map $i\colon C\to \{0,1\}$ determined by $z\in C_{i(z)}$. Hence, if one
wants to prove the non-existence of almost invariant, non-trivial
colorings, it suffices to consider colorings where two colors
are used.

If $S\subset G$ is a set of generators for $G$, a coloring is almost
invariant if and only if for each $g\in S$ we have $c(x) = c(gx)$ for
all but finitely many $x\in X$. This observation is of course
particularly useful when $G$ is finitely generated, which will be the
case in this paper. Hence both $G$ and $X$ are countable. Also, it is
easy to see that any almost invariant coloring of $X$ can at most use
finitely many colors: Assume WLOG that $c\colon X\to C$ is surjective,
and for $z\in C$ let $X_z = c\inv(z)$; then $X = \bigsqcup_{z\in C}
X_z$ is the partition of $X$ associated to $c$.  Next, choose some
finite set of generators $g_1, \dots, g_k $ of $G$.  The almost
invariance of the coloring implies that each $g_i$ acts as a
permutation of all but finitely many $X_z$, hence $G$ acts as a
permutation on all but finitely many of the subsets. If the partition
consists of infinitely many subsets, this contradicts the assumption
that $G$ acts transitively on $X$.

Let $\Sigma$ be an oriented connected surface of genus $g$ with $r$
boundary components, where $g\geq 2$ and $r\geq 0$, or $g=1$ and
$r=0$. Let $\Gamma$ be the mapping class group of $\Sigma$, the group
of orientation preserving diffeomorphisms of $\Sigma$ fixing the
boundary (if any) pointwise, modulo the group of diffeomorphisms
isotopic (through isotopies fixing the boundary point-wise) to the
identity.  Furthermore, let $D_0$ be a non-empty \emph{multicurve} on
$\Sigma$, the isotopy class of a collection of disjoint circles in
$\Sigma$ which are not trivial nor parallel to a boundary component
(components of the multicurve are allowed to be parallel to each
other). Let $X$ denote the mapping class group orbit of $D_0$. 
%Our main results are these.

The main theorem of this paper is
\begin{theorem}
  \label{thm:2}
  When $g\geq 2$, $r$ arbitrary, there are no non-trivial almost
  invariant colorings of $X$.
\end{theorem}
The proof of this is the contents of
Section~\ref{cha:proof-theor-refthm:2}.
We are also able to prove that the ``converse'' is true for the closed
torus:
\begin{theorem}
  \label{thm:1}
  When $g=1$, $r=0$, there exist almost invariant colorings of $X$
  using arbitrarily many colors.
\end{theorem}
In fact, we classify all such almost invariant colorings explicitly.

\chapter{Motivation}
\label{cha:motivation}

Before delving into the proofs of Theorems~\ref{thm:2}
and~\ref{thm:1}, let us explain the motivation for studying the
question of existence of such almost invariant partitions. As in the
introduction, let $D_0$ be a multicurve on $\Sigma$, and let $X$ be
the mapping class group orbit of $D_0$. Let $M = \setC X$ denote the
complex vector spanned by $X$ (the set of finite formal $\setC$-linear
combinations of elements of $X$), and let $\hat M$ denote the
algebraic dual of $M$, which we may think of as the space of
\emph{all} formal linear combinations of elements of $X$. Both $M$ and
$\hat M$ become modules over $\Gamma$ by extending the action
$\setC$-linearly, and there is a $\Gamma$-equivariant inclusion
$\iota\colon M\hookrightarrow \hat M$.

In \cite{0710.2203}, we gave jointly with Andersen an algorithm for
computing the first cohomology group $H^1(\Gamma, \hat M)$ for any
multicurve $D_0$, and proved that for any surface $\Sigma$, there
exists some multicurve (in that paper called a \emph{BFK-diagram})
such that this cohomology group is non-trivial. We also gave an
explicit example showing that in the case of a closed torus, the
cohomology group $H^1(\Gamma, M)$ is non-trivial.

The short exact sequence
\begin{align}
  \label{eq:3}
  \cxymatrix{{0 \ar[r] & M \ar[r]^{\iota} & \hat M \ar[r] & \hat M/M
      \ar[r] & 0}}
\end{align}
of $\Gamma$-modules induces a long exact sequence of cohomology groups
\begin{align}
  \label{eq:4}
  \cxymatrix{{0 \ar[r] & H^0(\Gamma, M) \ar[r]^{\iota_*} & H^0(\Gamma,
      \hat M) \ar[r] & H^0(\Gamma, \hat M/M) \\
    \ar[r] & H^1(\Gamma, M) \ar[r]^{\iota_*} & H^1(\Gamma, \hat M)
    \ar[r] & \cdots}}
\end{align}
Since $X$ is in general infinite (except when $D_0$ is the empty
multicurve), it is easy to see that $H^0(\Gamma, M) = 0$ and
$H^0(\Gamma, \hat M) = \setC$, since no finite linear combination of
elements of $X$ is invariant under $\Gamma$, while the constant linear
combinations are invariant elements of $\hat M$. 

Now, what is an invariant element of the quotient module $\hat M/M$?
It is represented by an element $\hat m$ of $\hat M$ such that for
each $\gamma\in \Gamma$, we have $\gamma \hat m = \hat m$ in $\hat
M/M$, or in other words $\gamma \hat m - \hat m \in M$ for each
$\gamma\in\Gamma$ (this is by the way exactly how the connecting
homomorphism in \eqref{eq:4} above is defined). Hence, thinking of an
element of $\hat M$ as a coloring of $X$ by complex numbers, we see
that an invariant element of $\hat M/M$ is represented by an almost
invariant $\setC$-coloring of $X$. In terms of the exact sequence
\eqref{eq:4} above, Theorem~\ref{thm:2} implies that $H^0(\Gamma, \hat
M/M) = \setC$ and hence that the second $\iota_*$ is injective
whenever $g\geq 2$. So by computing the image of this map we obtain a
computation of $H^1(\Gamma, M)$, and this is done in \cite{0802.4372},
using the description of $H^1(\Gamma, \hat M)$ given in
\cite{0710.2203} and methods similar to those applied in the present
paper.

The study of the cohomology groups $H^1(\Gamma, M)$ is in turn
motivated by the fact that the complex vector space spanned by the set
of all (isotopy classes of) multicurves is isomorphic to the space
$\mathcal{O} = \mathcal{O}(\mathcal{M}_{\SL_2(\setC)})$ of algebraic
functions on the moduli space of flat $\SL_2(\setC)$ connections over
$\Sigma$. As a $\Gamma$-module, $\mathcal{O}$ splits into a direct sum
\begin{align}
  \label{eq:5}
  \mathcal{O} = \bigoplus_{D} M_D
\end{align}
where $M_D = \setC(\Gamma D)$ is the complex vector space spanned by
the $\Gamma$-orbit through $D$ and the sum is taken over a set of
representatives of the $\Gamma$-orbits. The splitting \eqref{eq:5}
then induces a splitting in cohomology
\begin{align}
  \label{eq:6}
  H^1(\Gamma, \mathcal{O}) = \bigoplus_D H^1(\Gamma, M_D).
\end{align}
Hence, by combining the results of the present paper with those of
\cite{0710.2203} we obtain in \cite{0802.4372} a complete calculation
of the left-hand side of \eqref{eq:6}.

\chapter{Proof of Theorem~\ref{thm:2}}
\label{cha:proof-theor-refthm:2}

\section{Useful facts}
\label{sec:useful-facts}

We are going to need a couple of facts regarding the mapping class
group and its action on the set of multicurves. First of all, the
mapping class group is generated by Dehn twists, and moreover there
exists a finite set of % non-separating
curves such that the Dehn twists on these curves generate $\Gamma$.
Furthermore, one may choose these curves so that any pair of them
intersect in at most two points (see \cite{MR1851559}). Dehn twists on
disjoint curves commute. When $g\geq 2$, a twist on a separating curve
can be written as a product of twists on non-separating curves. Hence
in this case the mapping class group is generated by a finite set of
twists in non-separating curves (though we may not necessarily choose
this set so that each pair of curves intersect in at most two points).

There is simple way to parametrize the set of all multicurves which
was found by Dehn. For details, we refer to~\cite{MR1144770}.
Essentially one cuts the surface into pairs of pants using $3g+r-3$
simple closed curves $\gamma_k$, and for each pair of pants one
chooses a set of three disjoint arcs connecting the three pairs of
boundary components.  Then for each pants curve $\gamma_k$ one records
the geometric intersection number $m_k(D) = i(\gamma_k, D)$ (which is
a non-negative integer) and a ``twisting number'' $t_k(D)$, which can
be any integer. This defines a $6g+2r-6$-tuple of integers $(m_1(D),
t_1(D), \ldots, m_{3g+r-3}(D), t_{3g+r-3}(D))$ (satisfying certain
conditions), and, conversely, from any such tuple satisfying these
conditions one may construct a multicurve.

The important fact is that in this parametrization, the action of the
twist in the curve $\gamma_k$ on a multicurve $D$ is given by
\begin{align}
  \label{eq:2}
  t_k(\tau_{\gamma_k}^{\pm 1} D) = t_k(D) \pm m_k(D),
\end{align}
all other coordinates being unchanged. The formula \eqref{eq:2} is
intuitive in the sense that it says that for each time $D$ intersects
$\gamma_k$ essentially, the action of $\tau_{\gamma_k}$ on $D$ adds
$1$ to the twisting number of $D$ with respect to $\gamma_k$. This can
be used to prove a number of important facts.

\begin{lemma}
  \label{lem:8}
  Let $\gamma$ be a simple closed curve and $D$ a multicurve. Then the
  following are equivalent:
  \begin{enumerate}[(1)]\firmlist
  \item The twist $\tau_\gamma$ acts trivially on $D$.
  \item The twist $\tau_\gamma$ acts trivially on each component of $D$.
  \item The geometric intersection number between $\gamma$ and $D$ is zero.
  \item One may realize $\gamma$ and $D$ disjointly.
  \end{enumerate}
  Conversely, if $\tau_\gamma$ acts non-trivially on $D$, all the
  multicurves $\tau_\gamma^n D$, $n\in\setZ$, are distinct.
\end{lemma}
\begin{proof}
  All of the above assertions can be proved from~\eqref{eq:2} by
  letting $\gamma$ be part of a pants decomposition of the surface.
  This is clearly possible if $\gamma$ is non-separating, while if
  $\gamma$ is separating, observe that both connected components
  resulting from cutting along $\gamma$ must have negative Euler
  characteristic (otherwise $\gamma$ would be trivial or parallel to a
  boundary component, in which case the twist on $\gamma$ clearly acts
  trivially on $D$).
\end{proof}
To find a twist acting non-trivially on a multicurve, we need only
find a curve which has positive geometric intersection number with
the multicurve. This is possible if and only if the multicurve has a
component which is not parallel to a boundary component of $\Sigma$.

On a surface with negative Euler characteristic, there exist complete
hyperbolic metrics of constant negative curvature. Within each free
homotopy class of simple closed curves, there is a unique geodesic
representative with respect to such a metric. If $a$ and $b$ are the
geodesic representatives of distinct homotopy classes $\alpha$,
$\beta$, then $a$ and $b$ realizes the geometric intersection number
between $\alpha$ and $\beta$, ie. $\# a\cap b = i(\alpha, \beta)$.

% A Dehn twist acts trivially on a multicurve if and only if it acts
% trivially on each component of the multicurve. In other words, a Dehn
% twist can not permute the components of a multicurve. In fact, given
% simple closed curves $gamma$ and $\alpha$, let $i = i(\gamma, \alpha)$
% denote the geometric intersection number of $\gamma$ and $\alpha$.
% Then if one chooses some orientations of the two curves we have in
% homology that $\tau_\gamma \alpha = \alpha + i\gamma$. This also
% proves another important fact: If $\tau_\gamma$ acts non-trivially on
% a multicurve $D$, then $\tau_\gamma^n D \not= \tau_\gamma^m D$ for
% $n\not= m$.

\section{Interesting pairs}
\label{sec:interesting-pairs}

We will assume that the elements of $X$ have been colored red and
blue, and then prove that one of these colors has only been used a
finite number of times. To this end, an \emph{interesting pair} is a
pair $(\tau_\gamma, D)$ where $\tau_\gamma$ is a Dehn twist in a
% non-separating 
curve $\gamma$ and $D\in X$ is a multicurve such that
$\tau_\gamma D \not= D$ (equivalently, $i(\gamma, D) > 0$). Since
$\tau_\gamma$ changes the color of only finitely many diagrams, the
diagrams $\tau_\gamma^n D$ all have the same color for all
sufficiently large values of $n$. This color is called the
\emph{future} of the interesting pair $(\tau_\gamma, D)$, denoted
$\fut(\tau_\gamma, D)$. Similarly, we may consider the \emph{past}
$\pas(\tau_\gamma, D)$ of an interesting pair; the common color of all
diagrams $\tau_\gamma^{-n} D$ for sufficiently large $n$. We will also
need to consider pairs of the form $(\tau_\gamma\inv, D)$; the same
definition of future and past applies to these, and clearly
$\fut(\tau_\gamma^{\pm1}, D) = \pas(\tau_\gamma^{\mp1}, D)$.
\begin{lemma}
  \label{lem:1}
  For any interesting pair $(\tau_\alpha, D)$, we have
  \begin{align}
    \label{eq:7}
    \pas(\tau_\alpha\inv, D) = \fut(\tau_\alpha, D) = \pas(\tau_\alpha, D) = \fut(\tau_\alpha\inv, D) 
  \end{align}
\end{lemma}
\begin{proof}
  It suffices to prove the middle identity. We may find a
  non-separating simple closed curve $\beta$ different and disjoint
  from $\alpha$ such that $(\tau_\beta, D)$ is also interesting. To
  see this, let $\delta$ be a component of $D$ for which $\tau_\alpha
  \delta \not= \delta$, and assume that $\alpha$ and $\delta$ are
  represented by geodesics with respect to some choice of hyperbolic
  metric. Cutting $\Sigma$ along $\alpha$ then yields a (possibly
  non-connected) surface with geodesic boundary, in which $\delta$ is
  a number of properly embedded hyperbolic arcs. At least one of the
  connected components of the cut surface has genus at least $1$, so
  in this component we may find a closed geodesic $\beta$, not
  parallel to a boundary component, intersecting one of the
  $\delta$-arcs. In the original surface, $\beta$ is still a geodesic
  intersecting the geodesic $\delta$; hence $\tau_\beta \delta \not=
  \delta$ and $(\tau_\beta, D)$ is interesting.

  Next, since $\tau_\alpha$ and $\tau_\beta$ commute, we see that
  $\tau_\alpha^n \tau^m D$ is an $\setZ\times\setZ$-indexed family of
  distinct multicurves. By assumption, both $\tau_\alpha$ and
  $\tau_\beta$ change the color of finitely many multicurves. Hence,
  outside some bounded region in $\setZ\times\setZ$, moving from one
  diagram to a neighbour does not change the color, and since we can
  connect the future of $(\tau_\alpha, D)$ to its past using such
  moves, the claim follows.
\end{proof}
From now on, we will only consider the future.
\begin{lemma}
  \label{lem:7}
  Assume that $\alpha$ and $\beta$ are simple closed curves with
  $i(\alpha,\beta)\leq 1$, and that $D$ is a multicurve such that
  $(\tau_\alpha, D)$, $(\tau_\beta, D)$ are interesting pairs. Then
  $\fut(\tau_\alpha, D) = \fut(\tau_\beta, D)$.
\end{lemma}
\begin{proof}
  If $i(\alpha, \beta) = 0$ the result follows from the proof of
  Lemma~\ref{lem:1}.

  Now assume $i(\alpha, \beta) = 1$. Then $\alpha \cup \beta$ is
  contained in a subsurface $\Sigma'$ of genus $1$ with one boundary
  component $\gamma$. If $D$ can not be isotoped to be contained
  entirely in $\Sigma'$, either some component of $D$ intersects
  $\gamma$ essentially, or some component of $D$ lives in the
  complement of $\Sigma'$. In the former case, it is clear that
  $(\tau_\gamma, D)$ is interesting, so the $i=0$ case implies
  $\fut(\tau_\alpha, D) = \fut(\tau_\gamma, D) = \fut(\tau_\beta, D)$.
  In the latter case, use the fact that the complement of $\Sigma'$
  has genus at least $1$ to find a simple closed curve intersecting
  $D$ essentially.

  Otherwise, $D$ lives entirely in $\Sigma'$. Let $D_\alpha$ denote
  any component of $D$ on which $\tau_\alpha$ acts non-trivially. Then
  $D_\alpha$ is a simple closed curve in a torus with one boundary
  component. Since $D_\alpha$ is not a parallel copy of the boundary
  component, it must be a non-separating curve not parallel to
  $\alpha$. Hence, thinking of $\alpha$ as a $(1,0)$-torus knot and
  $\beta$ as a $(0,1)$-torus knot, we conclude that $D_\alpha$ is a
  $(p, q)$-torus knot with $(p,q)\not= (1,0)$. But then any other
  component of $D$ is forced to be either parallel to the boundary
  component of $\Sigma'$ or to $D_\alpha$. The only way that
  $\tau_\beta$ can act on some component of $D$ is then that
  $\tau_\beta$ acts on $D_\alpha$; hence also $(p, q) \not= (0,1)$.

  Consider the schematic picture of $\Sigma'$ on
  Figure~\vref{fig:sigma-prime}, where the boundary component is the
  circle in the center and $\alpha$ and $\beta$ are the sides of the
  square.

  \setlength{\captionmargin}{2cm}
  \begin{figure}[htb]
%    \begin{adjustwidth}{2cm}{2cm}
      \centering
      \includegraphics{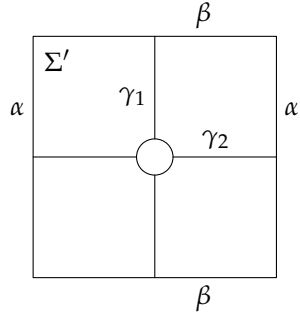}
      \caption{A torus with one boundary component.}
      \label{fig:sigma-prime}
%    \end{adjustwidth}
  \end{figure}

  We construct two disjoint simple closed %non-separating
  curves $\gamma_1, \gamma_2$ as follows: Draw two essential, disjoint
  arcs in $\Sigma'$ with the endpoints on the boundary component, and
  use the fact that the complement of $\Sigma'$ has genus at least $1$
  to close them up in such a way that they are disjoint and not
  homotopic to a curve contained in $\Sigma'$. 
  % By construction,
  % $\gamma_1$ intersects $\beta$ in exactly one point, and similarly
  % $\gamma_2$ intersects $\alpha$ in exactly one point, so they are
  % both non-separating.
  By the above description of $D_\alpha$, $(\tau_{\gamma_n}, D)$
  are both interesting pairs. Now the $i=0$ case implies that
  \begin{align*}
%    \label{eq:1}
    \fut(\tau_\alpha, D) = \fut(\tau_{\gamma_1}, D) =
    \fut(\tau_{\gamma_2}, D) = \fut(\tau_\beta, D).
  \end{align*}
\end{proof}
The next proposition extends the above lemma to $i(\alpha, \beta) \leq
2$, but its proof is rather technical. Also, as explained in the
comments following the proof, it is in fact not needed when one is
only interested in surfaces with at most one boundary component.
\begin{proposition}
  \label{prop:4}
  Assume that $\alpha$ and $\beta$ are simple closed curves with
  $i(\alpha,\beta) = 2$, and that $D$ is a multicurve such that
  $(\tau_\alpha, D)$ and $(\tau_\beta, D)$ are interesting. Then
  $\fut(\tau_\alpha, D) = \fut(\tau_\beta, D)$.
\end{proposition}
\begin{proof}
  Let $N$ be a regular neighbourhood of $\alpha\cup \beta$. We
  distinguish between four cases.
  \begin{enumerate}[(1)]\firmlist
  \item At least one of $\alpha$ and $\beta$ is non-separating in $N$.
  \item Both $\alpha$ and $\beta$ are separating in $N$, but
    non-separating in $\Sigma$.
  \item Both $\alpha$ and $\beta$ are separating in $N$, but one is
    non-separating in $\Sigma$.
  \item Both $\alpha$ and $\beta$ are separating in $\Sigma$.
  \end{enumerate}
  
  In case (1), assume WLOG that $\alpha$ is non-separating. This means
  that when cutting $N$ along $\alpha$, there is at least one arc $b$
  of $\beta$ connecting the two sides of $\alpha$. Now construct two
  curves $\gamma_1$, $\gamma_2$ as follows: Make two parallel copies
  of $b$ and close them up using arcs going in opposite directions
  along $\alpha$.  Applying small isotopies in a tubular neighbourhood
  of $\alpha$ we obtain a situation as depicted in
  Figure~\vref{fig:one-non-separating}. We observe that each
  $\gamma_n$ intersects $\alpha$ in exactly one point, and also they
  intersect each other in exactly one point $p$. Furthermore, since
  $i(\alpha, \beta) = 2$, the arc $b$ does not start and end at the
  same point of $\alpha$, so we have $i(\gamma_n, \beta) = 1$ for $n=1,2$.

  \setlength{\captionmargin}{2cm}
  \begin{figure}[htb]
      \centering
      \includegraphics{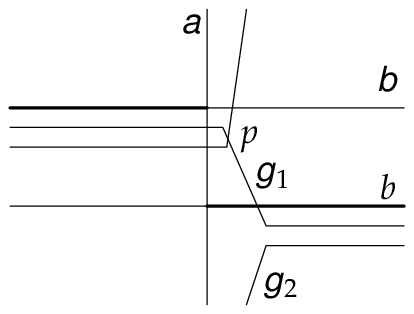}
      \caption{When $\alpha$ is non-separating in $N$, the two sides of
        $\alpha$ are connected by an arc of $\beta$.}
      \label{fig:one-non-separating}
  \end{figure}
  Now orient $\gamma_1$ and $\gamma_2$ oppositely along $b$. Then
  Goldman's bracket (see \cite{MR846929}) of $\gamma_1$ and $\gamma_2$
  is plus or minus some oriented version $\vec\alpha$ of $\alpha$.
  Now let $D_\alpha$ be some component of $D$ on which $\tau_\alpha$
  acts non-trivially. We claim that at least one of $\gamma_1$ and
  $\gamma_2$ intersects $D_\alpha$ essentially.  If this were not the
  case, choose geodesic representatives $\gamma_1'$, $\gamma_2'$ and
  $D_\alpha'$ of the three curves. Then $\gamma_i'$ is disjoint from
  $D_\alpha'$, and necessarily $\gamma_1'$ and $\gamma_2'$ intersect
  transversally in a single point $p'$. But then
  $(\gamma_1'\gamma_2')_{p'} \in \pi_1(\Sigma, p')$ is a
  representative of the free homotopy class of $\vec \alpha$ which
  does not intersect $D_\alpha'$, implying that $i(D_\alpha, \alpha) =
  0$, which contradicts the choice of $D_\alpha$.
  % If this were not the case, then by the Jacobi identity
  % \begin{align*}
  %   0 = [[\gamma_1, \gamma_2], D_\alpha] + [[\gamma_2, D_\alpha],
  %   \gamma_1] + [[D_\alpha, \gamma_1], \gamma_2] = [\pm \vec\alpha, D_\alpha]
  % \end{align*}
  % which contradicts the fact that the Goldman bracket counts the
  % number of essential intersection points between two simple curves
  % when one has a simple representative (see \cite{MR2041630} and
  % \cite{07062439}).
  So one of the pairs $(\tau_{\gamma_n}, D)$ is interesting, and by
  Lemma~\ref{lem:7} we have
  \begin{align*}
    \fut(\tau_\alpha, D) = \fut(\tau_{\gamma_n}, D) = \fut(\tau_\beta,
    D).
  \end{align*}
  This ends case (1).

  In cases (2)--(4), notice that $N$ is necessarily a sphere
  with four holes, and $\alpha$ and $\beta$ divide $N$ into two pairs
  of pants in two different ways. Denote the boundary components of
  $N$ by $\gamma_i$, $i=0,1,2,3$, such that $\gamma_1, \gamma_2$ are
  on one side of $\alpha$ and $\gamma_3, \gamma_4$ on the other, and
  such that $\gamma_2, \gamma_3$ are on one side of $\beta$ and
  $\gamma_4, \gamma_1$ on the other. Schematically we have
  Figure~\vref{fig:case-two--four}.

  Throughout the rest of the proof, we assume that $\alpha$, $\beta$,
  $\gamma_i$, $i=0,1,2,3$, denote geodesic representatives for their
  isotopy classes. Also, we let $\delta$ be the geodesic
  representative of some component of $D$ on which $\tau_\alpha$ acts
  non-trivially.  If $\delta$ does not live entirely in $N$, a twist
  in one of the boundary components acts non-trivially on $\delta$,
  and since this boundary component is disjoint from $\alpha$ and
  $\beta$ we are done by Lemma~\ref{lem:7}.  Otherwise, $\delta$ is a
  separating curve in $N$ which is not parallel to a boundary
  component. Clearly $\delta$ can not be parallel to $\beta$, since in
  that case $D$ could not consist of any component on which
  $\tau_\beta$ acts non-trivially. Hence $\delta$ is different from
  both $\alpha$ and $\beta$.

  In case (2), it is not hard to see that at least one of the
  ``opposite'' pairs $\gamma_1, \gamma_3$ and $\gamma_2, \gamma_4$ can
  be connected by an arc in the complement of $N$. Take two parallel
  copies of this arc, and close them up by arcs intersecting each
  other, $\alpha$ and $\beta$ exactly once as in
  Figure~\vref{fig:case-two} (the two connecting arcs are related by a
  twist in $\alpha$. We may then argue exactly as in case (1)
  to see that the twist in at least one of these simple closed curves
  acts non-trivially on the multicurve in question.
  \setlength{\captionmargin}{0pt}
  \begin{figure}[htb]
    \begin{adjustwidth}{1cm}{1cm}
      \centering
      \subfloat[\label{fig:case-two--four} $N$ is a sphere with four
      holes.]{\includegraphics[trim = -15 0 -15 0,clip]{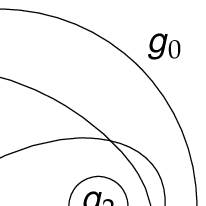}}
      \hfill
      \subfloat[\label{fig:case-two} In case (2), two opposite boundary
      components are connected in the complement of
      $N$.]{\includegraphics[trim = -15 0 -15 0,clip]{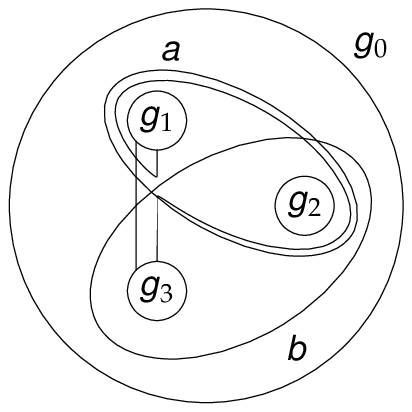}}
      
      \subfloat[\label{fig:case-three} In case (3), two ``neighbouring''
      boundary components are connected in the complement of
      $N$.]{\includegraphics[trim = -15 0 -15 0,clip]{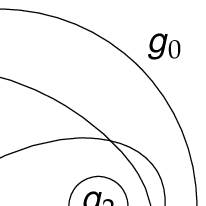}}
      \hfill
      \subfloat[\label{fig:case-four} In case (4), there exists an
      essential arc in the complement of $N$ starting and ending at the
      same boundary component.]{\includegraphics[trim = -15 0 -15
        0,clip]{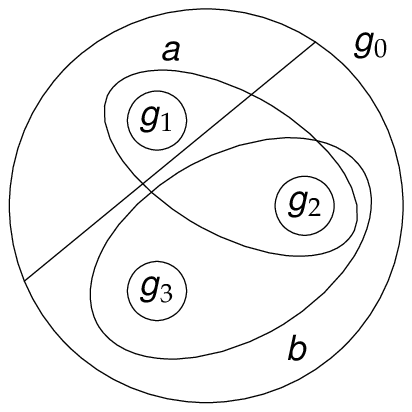}}
    \end{adjustwidth}
    \setlength{\captionmargin}{1cm}
    \caption{There are four different topological cases when two
        curves intersect in two points.}
  \end{figure}

  In case (3), assume WLOG that $\beta$ is separating and $\alpha$ is
  non-separating. This means that it is impossible to connect any of
  $\gamma_0$ and $\gamma_1$ to any of $\gamma_2$ and $\gamma_3$ in the
  complement of $N$. But then, since $\alpha$ is non-separating, one
  may connect either $\gamma_0$ to $\gamma_1$ or $\gamma_2$ to
  $\gamma_3$ in the complement of $N$. Assume WLOG that the latter is
  the case, and construct a simple closed curve $\gamma$ disjoint from
  $\beta$ intersecting $\gamma_2$, $\alpha$ and $\gamma_3$ exactly
  once each by composing the arc in the complement of $N$ with an arc
  in $N$, as in Figure~\vref{fig:case-three}. Observe that the
  geodesic representative of $\gamma$ necessarily intersects
  $\gamma_2$, $\alpha$ and $\gamma_3$ exactly once and is disjoint
  from $\beta$, so this representative contains a subarc in $N$
  starting at $\gamma_2$ and ending at $\gamma_3$. We now claim that
  this arc intersects $\delta$ (recall that $\delta$ has been chosen
  to be a geodesic). Assume the contrary. Then $\delta$ is a simple
  closed curve in the surface obtained by cutting $N$ along this arc,
  which is a pair of pants. The ``legs'' are $\gamma_0$ and
  $\gamma_1$, whereas the ``waist'' is composed of four segments; two
  copies of the connecting arc and the remaining boundary components
  (cut open). Since $\delta$ is simple, it is parallel to one of the
  boundary components of the pair of pants. But $\delta$ is certainly
  not parallel to any of the original boundary components, nor is it
  parallel to the ``waist'', since the latter is parallel to $\beta$.
  This contradiction implies that $(\tau_\gamma, D)$ is an interesting
  pair, and since $\gamma$ is disjoint from $\beta$ and intersects
  $\alpha$ in a single point, Lemma~\ref{lem:7} yields the desired
  result,
   \begin{align*}
    \fut(\tau_\alpha, D) = \fut(\tau_\gamma, D) = \fut(\tau_\beta, D).
  \end{align*}

  Finally, in case (4), none of the four boundary components of $N$
  can be connected in the complement of $N$. This means that at least
  one of the connected components of $\Sigma - N$ must have positive
  genus. Assume WLOG that the component $\Sigma_0$ bounded by
  $\gamma_0$ has positive genus. Now take some non-separating,
  essential arc in $\Sigma_0$ with its endpoints on $\gamma_0$ and
  compose it with some essential arc in $N$ disjoint from $\beta$ and
  intersecting $\alpha$ in exactly two points (cf.
  Figure~\ref{fig:case-four}) to obtain a non-separating curve
  $\gamma$ in $\Sigma$. We claim that $\tau_\gamma$ acts non-trivially
  on $\delta$, ie. that the arc in $N$ intersects $\delta$
  essentially. To see this, we argue as in case (3) above. Observe
  that $\gamma$ has geometric intersection number $2$ with $\alpha$
  and $\gamma_0$. Hence, the geodesic representative of $\gamma$
  intersects $\alpha$ and $\gamma_0$ exactly twice, so this geodesic
  contains a subarc in $N$ looking as the one depicted in
  Figure~\ref{fig:case-four}. We claim that this arc intersects
  $\delta$. If this were not the case, we may cut $N$ along this arc
  to obtain a cylinder (bounded by one of the original boundary
  components and a curve coming from the cut) and a pair of pants
  (bounded by two of the original boundary components and a curve from
  the cut), and $\delta$ lives completely in one of these.  Since
  $\delta$ is not parallel to any of the boundary components of $N$,
  we conclude that $\delta$ is parallel to the third boundary
  component of the pair of pants. But this third boundary component is
  clearly parallel to $\beta$, which contradicts the fact that $D$
  does not contain any component parallel to $\beta$. Hence
  $(\tau_\gamma, D)$ is interesting, and since $\gamma$ is
  non-separating and intersects $\alpha$ in two points, by case (3)
  and Lemma~\ref{lem:7} we have
   \begin{align*}
    \fut(\tau_\alpha, D) = \fut(\tau_\gamma, D) = \fut(\tau_\beta, D),
  \end{align*}
  which finishes the last case.
\end{proof}

Now we turn to the (finite) presentation of the mapping class group
given by Gervais in \cite{MR1851559}, where the generators are twists
in certain curves. A key property of this presentation is that any two
curves involved intersect each other in at most two points. 
% In the following, the $\eta$s will denote any of these curves.
It should be pointed out, however, that if one is only interested in
surfaces with at most one boundary component, a much earlier result by
Wajnryb \cite{MR719117} yields a presentation where each pair of
curves intersect in at most one point. In this case, one does not need
the rather technical Proposition~\ref{prop:4} above in the following
(simply replace all references to \cite{MR1851559} by \cite{MR719117}
and all occurences of ``at most two points'' by ``at most one
point'').
\begin{proposition}
  \label{prop:3}
  Let $S$ denote the set of curves from \cite{MR1851559} such that
  $\{\tau_\sigma \mid \sigma\in S\}$ generate $\Gamma$. Let
  $\alpha,\beta\in S$ be two of these curves, and let $D_1, D_2\in X$
  be multicurves such that $(\tau_\alpha, D_1)$ and $(\tau_\beta,
  D_2)$ are interesting. Then
  \begin{align*}
    \fut(\tau_\alpha, D_1) = \fut(\tau_\beta, D_2).
  \end{align*}
\end{proposition}
\begin{proof}
  We may find a sequence of curves $\eta_1, \eta_2, \dots, \eta_n\in
  S$ and exponents $\epsilon_i=\pm 1$ such that, writing $\tau_i =
  \tau_{\eta_i}^{\epsilon_i}$, $\tau_n\cdots\tau_2\tau_1 D_1 = D_2$.
  For each $1\leq i \leq n$ we may assume that $(\tau_i,
  \tau_{i-1}\cdots \tau_1 D_1)$ is interesting; otherwise we may
  simply omit the corresponding $\tau_i$. Now using alternately the
  fact that $\eta_i$ and $\eta_{i+1}$ intersect in at most two points
  and the obvious fact that $\fut(\tau_\gamma, D) = \fut(\tau_\gamma,
  \tau_\gamma D)$ for any interesting pair $(\tau_\gamma, D)$, we
  obtain a sequence of identities
  \begin{align*}
    \fut(\tau_1, D_1) &= \fut(\tau_1, \tau_1 D_1) = \fut(\tau_2,
    \tau_1 D_1)\\
    &= \fut(\tau_2, \tau_2\tau_1 D_1) = \fut(\tau_3, \tau_2\tau_1
    D_1)\\
    &\eqvdots\\
    &= \fut(\tau_{n-1}, \tau_{n-1}\cdots\tau_2\tau_1 D_1) =
    \fut(\tau_n, \tau_{n-1}\cdots\tau_2\tau_1 D_1)\\
    &= \fut(\tau_n, \tau_n\cdots\tau_2\tau_1 D_1) = \fut(\tau_n, D_2)
  \end{align*}
  which may be augmented by the identities $\fut(\tau_\alpha, D_1) =
  \fut(\tau_1, D_1)$ and $\fut(\tau_n, D_2) = \fut(\tau_\beta, D_2)$
  to obtain the desired result.
\end{proof}
\begin{lemma}
  \label{lem:6}
  Let $f\in \Gamma$ be any diffeomorphism, and $(\tau_\alpha, D)$ an
  interesting pair. Then $(\tau_{f(\alpha)}, fD)$ is also interesting
  and $\fut(\tau_\alpha, D) = \fut(\tau_{f(\alpha)}, fD)$.
\end{lemma}
\begin{proof}
  Recall that $f\circ \tau_\alpha \circ f\inv =
  \tau_{f(\alpha)}$. Hence $\tau_{f(\alpha)} (fD) = f(\tau_\alpha D)
  \not= fD$, so $(\tau_{f(\alpha)}, fD)$ is interesting. Also we have
  $\tau_{f(\alpha)}^n = f\circ \tau_\alpha^n \circ f\inv$, so
  $\tau_{f(\alpha)}^n (fD) = f(\tau_\alpha^n D)$. Since the different
  multicurves $\tau_\alpha^n D$ have the same color for all
  sufficiently large $n$, and since $f$ changes the color of only
  finitely many multicurves, the result follows.
\end{proof}

\begin{proposition}
  \label{prop:5}
  All interesting pairs $(\tau_\gamma, D)$ where $\gamma$ is a
  non-separating curve have the same future.
\end{proposition}
\begin{proof}
  Let $\tau_\alpha$ be a twist on a non-separating curve which is part
  of the generating set for $\Gamma$ from \cite{MR1851559}. Then
  Proposition~\ref{prop:3}, with $\alpha=\beta$, implies that the
  future is a property of $\tau_\alpha$ alone, and not of the
  particular multicurve on which $\tau_\alpha$ acts. If $\gamma$ is any
  non-separating curve, choose a diffeomorphism of $\Sigma$ carrying
  $\gamma$ to $\alpha$ and apply Lemma~\ref{lem:6}.
\end{proof}

Now we are ready to prove the main theorem.
\begin{proof}[of Theorem~\ref{thm:2}]
  Choose a finite set $\alpha_1, \dots, \alpha_N$ of
  non-separating curves such that the twists in these curves generate
  $\Gamma$ (we do not require that these intersect pair-wise in at
  most two points). To be concrete, assume that the common future (cf.
  Proposition~\ref{prop:5}) of all interesting pairs $(\tau_\gamma,
  D)$ with $\gamma$ non-separating is red. We must then prove that
  only finitely many multicurves are blue. Let $B\subset X$ be the set
  of blue multicurves. For each blue multicurve $D\in B$, choose a
  generator $\tau_{\alpha_k}$ such that $(\tau_{\alpha_k}, D)$ is
  interesting (this must be possible since the action is transitive
  and the $\tau_{\alpha_k}$ generate $\Gamma$). This defines a map
  $f\colon B\to \{1, 2, \dots, N\}$. We claim that for each $k\in\{1,
  \dots, N\}$, the pre-image $f\inv(k)$ is finite.

  To see this, for each $D\in f\inv(k)$ consider the
  ``$\tau_{\alpha_k}$-string through $D$'', ie. the set $s_k(D) =
  \{\tau_{\alpha_k}^n D \mid n\in\setZ\}$. Let $B_k$ be the union of
  the blue multicurves occuring in these strings, ie.
  \begin{align*}
    B_k = \bigcup_{D\in f\inv(k)} (s_k(D)\cap B),
  \end{align*}
  so that $f\inv(k) \subseteq B_k$. There are only finitely many blue
  multicurves in each string by Proposition~\ref{prop:5} and
  Lemma~\ref{lem:1}, and since $\tau_{\alpha_k}$ changes the color of
  at least one diagram in each string (since the strings contain both
  blue and red multicurves), there can be only finitely many strings
  by the almost invariance of the coloring. Hence, there are only
  finitely many blue multicurves.
\end{proof}

\chapter{The genus one case}
\label{cha:genus-one-case}

When $\Sigma$ is a closed torus, it is well-known that $\Gamma \cong
\SL_2(\setZ)$. A multicurve necessarily consists of some number of
parallel copies of the same non-separating simple closed curve, and
since $\Gamma$ acts identically on parallel curves, we may simply
assume that $D_0$ is a single non-separating simple closed curve, and
$X$ is the set of isotopy classes of such curves. Hence we may
identify $X$ with the set of unoriented torus knots, ie. the set $P$
of pairs $(p,q)$, $p,q\in\setZ$ and $\gcd(p,q) = 1$, where we identify
the pairs $(p,q)$ and $(-p, -q)$ (since the curves are not oriented).
The action of the mapping class group is then simply given by the
usual action of $\SL_2(\setZ)$ on pairs of relatively prime integers,
and the central element $-I$ acts trivially, so we are really dealing
with an action of $\PSL_2(\setZ)$.

As generators for $\SL_2(\setZ)$ we choose
$S = \left[
  \begin{smallmatrix}
    0 & 1\\ -1 & 0
  \end{smallmatrix}
\right]$
and
$R = \left[
  \begin{smallmatrix}
    0 & -1\\ 1 & 1
  \end{smallmatrix}
\right]$. Then $S^2 = R^3 = I$ in $\PSL_2(\setZ)$. Letting
\begin{align*}
  X_1 &= \{ (p,q) \mid p \geq 1, q\geq 0\}\\
  X_2 &= \{ (p,q) \mid q > -p \geq 0\}\\
  X_3 &= \{ (p,q) \mid -p \geq q > 0\}
\end{align*}
it is easy to see that $X_1 \cup X_2 \cup X_3 = X$, and one also
verifies that $S X_1 = X_2\cup X_3$, $R X_1 = X_2$, $R X_2 = X_3$.

\begin{figure}[htb]
  \centering
  \includegraphics{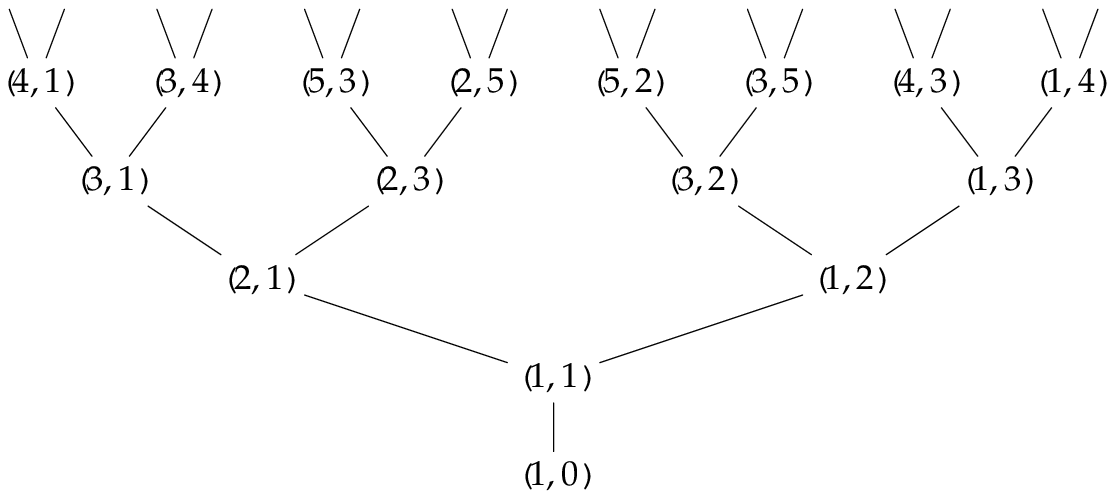}
  \caption{An infinite binary tree labelled by the points of $X_1$.}
  \label{fig:bin-tree}
\end{figure}

\begin{proposition}
  \label{prop:6}
  Any point $(p,q)\in X$ with $p,q>0$, can be reached from $(1,1)$ by
  applying a unique sequence of elements of $\SL_2(\setZ)$ of the form
  $S\inv R^k$, where $k$ is $1$ or $2$.
\end{proposition}
\begin{proof}
  For existence, we will use induction on $\max(p,q)$. For $\max(p,q)
  = 1$ we have $p=q=1$, in which case the claim is obvious (choose the
  empty sequence). If $\max(p,q) > 1$, $p$ and $q$ are different since
  $\gcd(p,q) = 1$. If $p>q$, put
  \begin{align*}
    (p',q') = R\inv S (p,q) = R\inv (-q, p) = (p-q, q)
  \end{align*}
  while if $q > p$, put
  \begin{align*}
    (p',q') = R^{-2} S (p,q) = R^{-2} (-q, p) = (p, q-p).
  \end{align*}
  In both cases, clearly $1\leq p',q'$ and $\max(p',q') < \max(p,q)$,
  so there exists $\gamma' = S\inv R^{k_{n-1}} \cdots S\inv R^{k_1}$
  with $\gamma'(1,1) = (p',q')$. Then $\gamma = S\inv R^{k_n} \gamma'$
  where $k_n = 1$ if $p>q$ and $k_n = 2$ if $p<q$ is an element of
  $\PSL_2(\setZ)$ of the desired form.

  To prove uniqueness, choose $(p,q)$ with $\max(p,q)$ minimal such
  that there are two different strings
  \begin{align*}
    \gamma_1 &= S\inv R^{k_n} S\inv R^{k_{n-1}} \cdots S\inv R^{k_1}\\
    \gamma_2 &= S\inv R^{\ell_m} S\inv R^{\ell_{m-1}} \cdots S\inv R^{\ell_1}
  \end{align*}
  satisfying $\gamma_i (1,1) = (p,q)$. Then $S(p,q)$ is a point in
  $X_2\cup X_3$ which is obtained by applying $R^{k_n}$ to some point
  of $X_1$ and also by applying $R^{\ell_m}$ to some (possibly other)
  point of $X_1$. But since $S(p,q)$ is an element of exactly one of
  $X_2 = R X_1$ and $X_3 = R^2 X_1$, this implies that $k_n =
  \ell_m$. Continuing this way, we only need to show that there is no
  non-trivial string
  \begin{align*}
    \gamma = S\inv R^{k_n} S\inv R^{k_{n-1}} \cdots S\inv R^{k_1}
  \end{align*}
  such that $\gamma(1,1) = (1,1)$. But this is trivial by observing
  that each element of the form $S\inv R^k$ strictly increases the
  $\max$-norm of any point $(p,q)$ with $p,q\geq 1$ (since $S\inv
  R(p,q) = (p+q, q)$ and $S\inv R^2(p,q) = (p,p+q)$).
\end{proof}
This proposition allows us to label the vertices of an infinite binary
tree $T$ as follows. The root is labelled by $(1,1)$, and all remaining
vertices are labelled according to the rule: If a vertex $v$ is
reached by going ``left'' from the immediate predecessor, the label of
$v$ is obtained by appying $S\inv R$ to the label of its predecessor;
otherwise the label is obtained by applying $S\inv R^2$ (see
Figure~\ref{fig:bin-tree}).

Now add a single vertex below the root and label this by $(1,0)$. This
gives, by Proposition~\ref{prop:6}, a $1$--$1$-correspondence between
the vertices of $T$ and the points in $X_1$, and from now on we shall
refer to a vertex and its label interchangeably.

By the level of a vertex of $T$ we mean its distance from $(1,1)$ (the
level of $(1,0)$ may be taken to be $-1$); there are $2^k$ vertices at
level $k$ for each $k\geq 0$, and also exactly $2^k$ vertices at level
$<k$.
\begin{proposition}
  \label{prop:7}
  For each $k\geq 0$, there is an almost invariant coloring of $X$
  using $2^k$ different colors.
\end{proposition}
\begin{proof}
  We start by coloring the subset $X_1$ by coloring the vertices of
  the tree. Assign different colors to the $2^k$ vertices at level
  $k$, and for each of these vertices assign the same color to all
  descendants. The remaining $2^k$ points of $X_1$ may be colored
  arbitrarily.

  To obtain a coloring of all of $X$, we insist that the coloring is
  completely invariant under $R$. This gives a well-defined coloring,
  since $X_1$ is a complete set of representatives of the $R$-orbits
  of $X$. In order to see that this coloring is almost invariant under
  $\PSL_2(\setZ)$, it suffices to check that the other generator $S$
  changes the color of only finitely many points of $X$. Since $S$ has
  order two in $\PSL_2(\setZ)$, $S$ changes the color of $p$ if and
  only if it changes the color of $Sp$. Hence we need only check that
  $S$ changes the color of finitely many elements of $X_1$. But for any
  vertex $v$ of $T$ of level $k+1$ or higher, applying $S$ to the
  label of $v$ yields by construction a point of $X_2$ or $X_3$ which
  has the same color as the predecessor of $v$; hence $S$ does not
  change the color of labels placed at level $k+1$ or higher, and thus
  $S$ changes the color of at most $2\cdot 2^{k+1}$ points of $X$.
\end{proof}
This finishes the proof of Theorem~\ref{thm:1}, but we also promised a
classification of all almost invariant colorings in this case. 
% We say that two almost invariant colorings of $X$ are equivalent if
% they differ on only finitely many elements of $X$. A coloring
% obtained from another coloring by identifying two or more colors is
% said to be a simplification of the latter.
\begin{proposition}
  \label{prop:8}
  Any almost invariant coloring of $X$ is equivalent to (a coloring
  which is a simplificaton of) a coloring of the form constructed in
  Proposition~\ref{prop:7}.
\end{proposition}
\begin{proof}
  Let $c\colon X\to C$ be some almost invariant coloring of $X$. Since
  $R$ changes the color of only finitely many points of $X$, it
  changes the color of only finitely many points $x_1, \dots, x_N$ of
  $X_1$. Now we change $c$ into an equivalent coloring $c'$ by putting
  $c'(Rx_i) = c'(R^2x_i) = c(x_i)$, and $c' = c$ otherwise. Then $c'$
  is by construction completely invariant under $R$. Now since $c'$ is
  almost invariant, there are only finitely many points of $X_1$ whose
  $c'$-color changes under $S$. Choose $K$ such that the color of any
  label placed at level $k>K$ is unchanged under $S$. This, together
  with the $R$-invariance of $c'$, implies that each label at level
  $K$ has the same color as any of its descendants, and hence $c'$ is
  (a simplification of) a coloring using $2^K$ different colors.
\end{proof}

\bibliographystyle{alphaurl}
\bibliography{../phd}

\noindent\emph{Acknowledgement:} The author would like to thank Jørgen
Ellegaard Andersen and Thomas Kragh for valuable discussions and
comments.

\end{document}